\documentclass[10pt]{article}

\usepackage{amsmath}
\usepackage{amssymb}
\usepackage{indentfirst}
\usepackage{graphics} 
\usepackage{color}

\makeatletter

\@addtoreset{equation}{section}
\makeatother
\def\<{\langle}
\def\>{\rangle}

\newtheorem{lem}{Lemma}[section]
\newtheorem{theo}{Theorem}[section]
\newtheorem{rem}{Remark}[section]
\newtheorem{pro}{Proposition}[section]

\makeatletter
   
   \@addtoreset{equation}{section}
\makeatother

\begin{document}
\title{\bf Local enegy decay for $2$-D wave equations\\ with variable coefficients}
\author{Ryo Ikehata\thanks{ikehatar@hiroshima-u.ac.jp} \\ {\small Department of Mathematics, Division of Educational Sciences}\\ {\small Graduate School of Humanities and Social Sciences} \\ {\small Hiroshima University} \\ {\small Higashi-Hiroshima 739-8524, Japan}}
\date{}
\maketitle
\begin{abstract}
This paper addresses the two-dimensional initial value problem in ${\bf R}^{2}$ for the wave equation with varying spatial coefficients in the main part. Assuming compactness in the support of the initial value, we report that the corresponding local energy decays to an order of magnitude of, for example, $O(t^{-1}\sqrt{\log t})$ after sufficiently large time. For the two-dimensional whole space case, it is crucial to establish the optimal $L^2$-estimate for the solution itself, skillfully avoiding the difficulty of not being able to use useful inequalities such as Hardy-type inequalities in higher dimensional case. We also consider cases where the variable coefficients are slightly generalized. These proofs are developed using the multiplier method.
\end{abstract}
\section{Introduction}
\footnote[0]{Keywords and Phrases: Wave equation, two dimension, bulk modulus, $L^{2}$-growth, local energy decay.}
\footnote[0]{2020 Mathematics Subject Classification. Primary 35L05; Secondary 35B40, 35C20.}

We are concerned with the Cauchy problem of the following $2$ dimensional wave equations with bulk modulus $K(x)$: 
\begin{equation}
u_{tt}(t,x) - \nabla\cdot(K(x)\nabla u(t,x)) = 0,\ \ \ (t,x)\in (0,\infty)\times {\bf R}^{2},\label{0eqn}
\end{equation}
\begin{equation}
u(0,x)= u_{0}(x),\ \ u_{t}(0,x)= u_{1}(x),\ \ \ x\in {\bf R}^{2},\label{0initial}
\end{equation}
where the initial data $u_{j}$ ($j = 0,1$) are taken from the test function space
\[u_{j} \in C_{0}^{\infty}({\bf R}^{2}) \quad (j = 0,1),\]
and satisfy the compact support condition such that 
\begin{equation}\label{I-1}
{\rm supp}\,u_{j} \subset B_{L}
\end{equation}
with some $L > 0$ for $j = 0,1$, where $B_{L} := \{x \in {\bf R}^{2}\,:\,\vert x\vert \leq L\}$, and (for simplicity)
\[K \in {\rm BC}^{1}({\bf R}^{2}).\]
All functions, solutions to equations, and coefficients appearing shall be real-valued.\\

{\bf Notation} {\small Throughout this paper, $\| \cdot\|$ stands for the $L^2({\bf R}^{2})$-norm, and we denote the usual $L^{p}$-norm of $u \in L^{p}({\bf R}^{2})$ by $\Vert u\Vert_{p}$ for $p \in [0,2)\cup (2,\infty]$. Sometimes the symbol $\Vert f\Vert_{L^{2}(\Omega)}$ is also used for $f \in L^{2}(\Omega)$ ($\Omega \subset {\bf R}^{2}$). $f \in {\rm BC}({\bf R}^{2}) \Leftrightarrow$ $f$ is continuous and bounded on ${\bf R}^{2}$. $f \in {\rm BC}^{1}({\bf R}^{2}) \Leftrightarrow$ $f \in {\rm BC}({\bf R}^{2}) $ and $\exists \partial f(x)/\partial x_{j} \in {\rm BC}({\bf R}^{2})$ ($j = 1,2$) for the function $f(x) = f(x_{1},x_{2})$.} \\

Under these conditions the Cauchy problem \eqref{0eqn}-\eqref{0initial} has a unique sufficiently smooth solution $u(t,x)$ that satisfies the energy conservation law:
 \[E_{u}(t) = E_{u}(0) \quad \forall t \in [0,\infty),\] 
 and the corresponding solution $u(t,x)$ has the finite speed of propagation property such that
\[{\rm supp}\,u(t,\cdot) \subset B_{L+k_{1}t},\]
where the total energy $E_{u}(t)$ for the solution $u(t,x)$ to problem \eqref{0eqn}-\eqref{0initial} is defined by
\begin{equation}
E_{u}(t) := \frac{1}{2}(\| u_t(t,\cdot)\|^2+\|\sqrt{K(\cdot)}\nabla u(t,\cdot)\|^2,
\end{equation}
and
$$k_{1} := \sup\{\sqrt{K(x)}\,:\,x \in {\bf R}^{2}\}.$$
For these facts, for example, one can refer to \cite{Ikawa}. For later use, for each $R > 0$ one also defines the localized energy of the solution $u(t,x)$ in the region $B_{R}$: 
\begin{equation}
E_{u,R}(t)=\frac{1}{2}\int_{B_{R}}\left(\vert u_t(t,x)\vert^2 + K(x)\vert\nabla u(t,x)\vert^2\right)dx.
\end{equation}

In this research, one imposes additionally three assumptions that \\ 
{\bf (K-1)}\,$K(x) \geq k_{m} > 0$ with some $k_{m} > 0$ for all $x \in {\bf R}^{2}$,\\
{\bf (K-2)}\, $x\cdot\nabla K(x) \leq 0$ for all $x \in {\bf R}^{2}$,\\ 
{\bf (K-3)}\,$K(x) = k_{0} (\geq k_{m})$ for all $x$ satisfying $\vert x\vert > r_{0}$ with some $r_{0} > 0$.\\

The topic of local energy decay in the wave equation
\begin{equation}\label{I-146}
u_{tt}(t,x)-\Delta u(t,x) = 0
\end{equation}
is a classical and fundamental problem, and it is also highly significant from the perspective of scattering theory. The first crucial point in this problem is Morawetz's derivation in \cite{Mora} of uniform decay in the order of $O(t^{-1})$ ($t \to \infty$) for a star-shaped exterior region in three-dimensional space (see also \cite{Lax}). The author in \cite{Mora} proposed the so-called Morawetz identity, as in Lemma \ref{2-1} with $K(x) = 1$, and solved the problem using the analysis on the corresponding Poisson equation combined with the multiplier method. Subsequently, the local energy decay in the exterior region of a general non-trapping obstacle was studied in \cite{Mora-3, Mora-2}. Furthermore, in \cite{Met} the recent integrated local energy estimate in a two-dimensional star-shaped exterior is also intriguing. It is noted that this poses difficulties when considering the problem in the $2D$ whole space without obstacles. Incidentally, note that these results essentially utilized the compactness of the support of the inital data. This paper's research analyzes based on the multiplier method, but for local energy decay utilizing spectral analysis, refer to, for example, \cite{Bu, Va, V} and the references therein. For further discussion on Morawetz's generalized problem independent of the initial value's compactness, see references such as \cite{IN, Mura, Z2}. Research in this direction appears to be still quite limited. In this context, in \cite{JHDE-ike}, the author recently examined the initial value problem for the free wave equation \eqref{I-146} in two-dimensional entire space. In \cite{JHDE-ike} the author derived that the corresponding local energy decays at the order of $O(t^{-1}\sqrt{\log t})$ (as $t \to \infty$) when the compactness of the initial value's support is not assumed. The two-dimensional entire space case is fundamentally difficult in the sense that, since the so-called Hardy-type inequalities used in \cite{IN, IS} do not hold (see e.g., \cite{MOW, L-I-2}), one must find alternative calculations. The essence of the proof lies in deriving the $L^2$ estimate formula for the solution itself. However, for the two-dimensional free wave case, as pointed out in \cite{JHDE-ike}, the fact that $\Vert u(t,\cdot)\Vert \sim \sqrt{\log t}$ ($t \to \infty$) holds is one factor contributing to its difficulty. For two-dimensional exterior problems, this is a significant difference since $\Vert u(t,\cdot)\Vert = O(1)$ ($t \to \infty$) holds fundamentally in the 2-D exterior domain case (see \cite{IM}). Fundamentally, these are analyses of mixed problems for wave equations outside obstacles. Conversely, research on local energy decay in the entire space without obstacles in low dimensions seems scarce. Moreover, analyses of general hyperbolic equations in “low dimensions” incorporating variable coefficients appear limited to references like \cite{CI, D, JHDE-ike, Sha} and the references therein. By the way, for reference, consult \cite{B, F, I-05, IS, Liu, Z} regarding local energy decay in hyperbolic equations with variable coefficients in exterior mixed problems. Among these, the result \cite{D} deriving the exponential decay in time of local energy for the initial value problem of the following equation in the whole one-dimensional space is particularly interesting:
\[\beta(x)w_{tt}(t,x) - (\alpha(x)w_{x}(t,x))_{x} = 0,\quad x \in {\bf R}.\]
Assuming compactness of the initial value support, such results are derived based on the spectral theory of self-adjoint operators. 

Furthermore, for waves with variable coefficients, the following equations (limited to two dimensions) derive logarithmic in time decay in \cite{Sha} and algebraic in time decay in \cite{CI}:
\[u_{tt}(t,x)-a(x)^{2}\Delta u(t,x) = 0, \quad x \in {\bf R}^{2}.\]
Note, however, that \cite{Sha} assumes compactness of the support of the initial data and high regularity on them, while \cite{CI} does not assume compactness of the support of the initial data but imposes strong conditions such that
\begin{equation}\label{2025}
\int_{{\bf R}^{2}}\frac{u_{1}(x)}{a(x)^{2}}dx = 0
\end{equation}
on the initial velocity $u_{1}(x)$.

Based on the above, this paper aims to determine the order of local energy decay for the two-dimensional whole-space initial value problem of the divergence form of wave equation with variable-coefficient \eqref{0eqn}, assuming only the compactness of the support of the initial data while omitting assumptions like \eqref{2025} regarding the zero-order moment of the initial velocity $u_{1}(x)$. Although the equation itself \eqref{0eqn} differs slightly, we also confirm that stronger assumptions like those imposed in \cite{CI} can be removed. However, as a trade-off, compactness must be assumed for the initial value support.

About the behavior as time goes to infinity one can state the following uniform decay of the localized energy.
\begin{theo}\label{theorem}\,Assume \eqref{I-1} on the initial data $u_{j} \in C_{0}^{\infty}({\bf R}^{2})$ {\rm (}$j = 0,1${\rm )} and the conditions {\bf (K-1)}, {\bf (K-2)} and {\bf (K-3)}. Let $R > r_{0}$ be an arbitrarily fixed real number. Then, for the solution $u(t,x)$ to problem \eqref{0eqn}-\eqref{0initial} it holds that
\[E_{u,R}(t) = O(t^{-1}\sqrt{\log t}) \quad (t \gg 1).\] 
\end{theo}
\begin{rem}{\rm Condition {\bf (K-1)} indicates that the main part of the equation satisfies the condition of ellipticity. Condition {\bf (K-2)} imposes a requirement that the bulk modulus $K(x)$ decreases monotonically in the radial direction. Roughly speaking, this represents a state where the medium's temperature is high near the origin and decreases with distance from the origin. In such cases, the local energy exactly decays uniformly with its rate $O(t^{-1}\sqrt{\log t})$ (as $t \to \infty$). Condition {\bf (K-3)} corresponds to the equation matching a free wave with constant coefficients in the spatial far field. This is a crucial condition, and it seems extremely difficult to obtain similar results if this far-field condition is omitted (see Remark \ref{K-2}).}
\end{rem}

This paper is organized as follows. In Section 2 one prepares the Morawetz's identity (\cite{Mora}), the Todorova-Yordanov identity (\cite{TY}), and the $L^2$-estimates for the solutions themselves in order to prove our main Theorem \ref{theorem}. Section 3 is devoted to the proof of Theorem \ref{theorem} based on preliminary results prepared in Section 2. Section 4 additionally examines the generalization of the condition {\bf (K-2)} imposed on $K(x)$. In Section 5 we consider the case when $K(x)$ is Lipschitz continuous on ${\bf R}^{2}$. \\

\section{Preliminary results}

We start with the so-called (generalized) Morawetz identity (see \cite{Mora} and \cite[Proposition 2.1]{I-05}). In Paper \cite{I-05}, the identity is obtained in the exterior region, but it holds in the entire space as well, with only the boundary integral term vanishing.
\begin{lem}\label{2-1} For the slution $u(t,x)$ to problem \eqref{0eqn}-\eqref{0initial} it holds that 
\[tE_{u}(t) + \frac{1}{2}\int_{{\bf R}^{2}}u_{t}(t,x)u(t,x)dx + \int_{{\bf R}^{2}}u_{t}(t,x)(x\cdot \nabla u(t,x))dx \]
\[= J_{0} + \frac{1}{2}\int_{0}^{t}\int_{{\bf R}^{2}}\left(x\cdot\nabla K(x)\right)\vert \nabla u(s,x)\vert^{2}dxds\]
for $t \geq 0$, where
\[J_{0} := \frac{1}{2}\int_{{\bf R}^{2}}u_{1}(x)u_{0}(x)dx + \int_{{\bf R}^{2}}u_{1}(x)(x\cdot\nabla u_{0}(x))dx.\]
\end{lem} 

One further prepares the following identity (see e.g., \cite[(2.8) at page 270]{IN} or \cite{IS}). This can be easily derived by the weighted energy method introduced originally by Todorova-Yordanov \cite{TY}.
For this we set pointwisely
\[e(t,x) := \frac{1}{2}\left(\vert u_{t}(t,x)\vert^{2} + K(x)\vert \nabla u(t,x)\vert^{2}\right).\]
\begin{lem}\label{2-2} Let $\eta \in C^{1}([0,\infty)\times {\bf R}^{2})$ satisfy $\eta_{t}(t,x) \ne 0$ for all $(t,x) \in [0,\infty)\times {\bf R}^{2}$. Then, for the slution $u(t,x)$ to problem \eqref{0eqn}-\eqref{0initial} it holds that
\[0 = \frac{\partial}{\partial t}\{\eta(t,x)e(t,x)\} -\nabla\cdot(\eta(t,x)u_{t}(t,x)K(x)\nabla u(t,x))\]
\[+ \frac{K(x)}{2(-\eta_{t}(t,x))}\left\vert \eta_{t}(t,x)\nabla u(t,x) - u_{t}(t,x)\nabla\eta(t,x) \right\vert^{2} \]
\[+ \frac{1}{2\eta_{t}(t,x)}\vert u_{t}(t,x)\vert^{2}\left(K(x)\vert\nabla\eta(t,x)\vert^{2} - \eta_{t}(t,x)^{2} \right) \]
for all $(t,x) \in [0,\infty)\times {\bf R}^{2}$.
\end{lem}

Now we prepare the following weight function $\psi:[0,\infty)\times{\bf R}^{2} \to {\bf R}$ as follows (see \cite{IS}).
\[ \psi(t,x) = \left\{
     \begin{array}{ll}
       \displaystyle{1 + \vert x\vert - \sqrt{k_{0}}t}&
              \qquad \vert x\vert \geq \sqrt{k_{0}}t, t \geq 0, \\[0.2cm]
       \displaystyle{(1 + \sqrt{k_{0}}t - \vert x\vert)^{-1}}& \qquad  \vert x\vert < \sqrt{k_{0}}t, t \geq 0.
\end{array} \right. 
\]      
Note that the function $\psi \in C^{1}([0,\infty)\times{\bf R}^{2})$ satisfies
\begin{equation}\label{2-3}
\psi_{t}(t,x) < 0,\quad \forall (t,x) \in [0,\infty)\times {\bf R}^{2}.
\end{equation}
\begin{equation}\label{2-3-1}
k_{0}\vert\nabla\psi(t,x)\vert^{2} - (\psi_{t}(t,x))^{2} = 0, \quad \forall (t,x) \in [0,\infty)\times {\bf R}^{2},
\end{equation}
\begin{equation}\label{2-3-2}
\psi(t,x) > 0,\quad \forall (t,x) \in [0,\infty)\times {\bf R}^{2}.
\end{equation}
Furthermore, for later use we prapare the auxiliary weight function $\phi \in C^{1}([0,\infty))$ as follows.
\[ \phi(t) = \left\{
     \begin{array}{ll}
       \displaystyle{1 + r_{0} - \sqrt{k_{0}}t}&
              \qquad r_{0} \geq \sqrt{k_{0}}t,\,\, t \geq 0, \\[0.2cm]
       \displaystyle{(1 + \sqrt{k_{0}}t - r_{0})^{-1}}& \qquad  r_{0} < \sqrt{k_{0}}t, \,\,t \geq 0.
\end{array} \right. 
\]  
We see that
\begin{equation}\label{2-3-3}
\phi(t) > 0,\quad \forall t \in [0,\infty),
\end{equation}
\begin{equation}\label{2-3-4}
\phi_{t}(t) < 0,\quad \forall t \in [0,\infty).
\end{equation}

One can prove the following weighted energy estimate.
\begin{lem}\label{2-6} There exists a constant $C_{r_{0}} > 0$ such that for the smooth solution $u(t,x)$ to problem \eqref{0eqn}-\eqref{0initial} it holds that
\[\int_{\vert x\vert \geq r_{0}}\psi(t,x)e(t,x)dx \leq C_{r_{0}}\int_{{\bf R}^{2}}(1+\vert x\vert)e(0,x)dx,\quad t \geq 0.\]
\end{lem} 
{\it Proof.}\,Let us prove the statement in order.\\

\underline{Weighted estimate on $B_{r_{0}}$}: By using \eqref{2-3-4}, let integrate both sides of the identity stated in Lemma \ref{2-2} with $\eta(t,x) := \phi(t)$ on $[0,t]\times B_{r_{0}}$. Then one can derive a series of inequalities.
\[0 \geq \int_{B_{r_{0}}}\phi(t)e(t,x)dx - \int_{B_{r_{0}}}\phi(0)e(0,x)dx - \int_{0}^{t}\int_{B_{r_{0}}}\nabla\cdot(\phi(s)u_{s}(s,x)K(x)\nabla u(s,x))dxds\]
\[+\int_{0}^{t}\int_{B_{r_{0}}}\frac{1}{2\phi_{s}(s)}\vert u_{s}(s,x)\vert^{2}\left(K(x)\vert\nabla\phi(s)\vert^{2} - \phi_{s}(s)^{2} \right)dxds \]
\[= \int_{B_{r_{0}}}\phi(t)e(t,x)dx - \int_{B_{r_{0}}}\phi(0)e(0,x)dx - \int_{0}^{t}\int_{B_{r_{0}}}\nabla\cdot(\phi(s)u_{s}(s,x)K(x)\nabla u(s,x))dxds\]
\[+\frac{1}{2}\int_{0}^{t}\int_{B_{r_{0}}}\vert u_{s}(s,x)\vert^{2}(- \phi_{s}(s))dxds \]
\begin{equation}\label{2-4}
\geq \int_{B_{r_{0}}}\phi(t)e(t,x)dx - \int_{B_{r_{0}}}\phi(0)e(0,x)dx - \int_{0}^{t}\int_{B_{r_{0}}}\nabla\cdot(\phi(s)u_{s}(s,x)K(x)\nabla u(s,x))dxds.
\end{equation}
Here, one has just used the fact that $\nabla\phi(t,x) = 0$ for $x \in B_{r_{0}}$ and $t \geq 0$.

\underline{Weighted estimate on $\vert x\vert \geq r_{0}$}:\,To get the weighted estimate in the region $\{\vert x\vert \geq r_{0}\}$, it should be noticed that the weight function $\psi(t,x)$ saisfies
\begin{equation}\label{K-1}
K(x)\vert\nabla\psi(t,x)\vert^{2} - \psi_{t}(t,x)^{2} = k_{0}\vert\nabla\psi(t,x)\vert^{2} - \psi_{t}(t,x)^{2} = 0
\end{equation}
for all $(t,x)$ satisfying $t \geq 0$ and $\vert x\vert > r_{0}$ (see the assumption {\bf (K-3)} and \eqref{2-3-1}). Then, integrating both sides of the identity stated in Lemma \ref{2-2} with $\eta(t,x) := \psi(t,x)$ on $[0,t]\times \{\vert x\vert \geq r_{0}\}$ and using \eqref{2-3} one has
\[0 \geq \int_{\vert x\vert \geq r_{0}}\psi(t,x)e(t,x)dx - \int_{\vert x\vert \geq r_{0}}\psi(0,x)e(0,x)dx\]
\[ - \int_{0}^{t}\int_{\vert x\vert \geq r_{0}}\nabla\cdot(\psi(s,x)u_{s}(s,x)K(x)\nabla u(s,x))dxds\]
\[+\int_{0}^{t}\int_{\vert x\vert \geq r_{0}}\frac{1}{2\psi_{s}(s,x)}\vert u_{s}(s,x)\vert^{2}\left(k_{0}\vert\nabla\psi(s,x)\vert^{2} - \psi_{s}(s,x)^{2} \right)dxds\]
\[= \int_{\vert x\vert \geq r_{0}}\psi(t,x)e(t,x)dx - \int_{\vert x\vert \geq r_{0}}\psi(0,x)e(0,x)dx \]
\begin{equation}\label{2-5}
- \int_{0}^{t}\int_{\vert x\vert \geq r_{0}}\nabla\cdot(\psi(s,x)u_{s}(s,x)K(x)\nabla u(s,x))dxds.
\end{equation}

Because of the divergence formula and the fact that $\phi(t) = \psi(t,x)$ for $\vert x\vert = r_{0}$ one notices that 
\[\int_{0}^{t}\int_{\vert x\vert \leq r_{0}}\nabla\cdot(\phi(s)u_{s}(s,x)K(x)\nabla u(s,x))dxds\]
\[ + \int_{0}^{t}\int_{\vert x\vert \geq r_{0}}\nabla\cdot(\psi(s,x)u_{s}(s,x)K(x)\nabla u(s,x))dxds = 0.\]
Thus, by adding both sides of \eqref{2-4} and \eqref{2-5}, and using \eqref{2-3-3} one can arrive at the desired estimate.
\hfill
$\Box$
\begin{rem}\label{K-2}
{\rm If a concrete form for the solution to the Eikonal equation \eqref{K-1} can be found, even just in the far-field region $\vert x\vert \geq r_{0}$, then the argument similar to the above proof becomes possible. Therefore, condition {\bf (K-3)} is absolutely indispensable at this point. }
\end{rem}
The following claim is a result specific to two dimensions. The fundamental idea was inspired by the recent paper \cite{RI}. The assumption \eqref{I-1} of compactness of the initial value support is used only in the proof of this proposition. The $L^2$ growth property of solutions to the two-dimensional wave equation is critical in a sense, and does not appear to be incorporated into frameworks such as \cite{St}.
\begin{pro}\label{2-7}
The smooth solution $u(t,x)$ to problem \eqref{0eqn}-\eqref{0initial} satisfies the following growth property:
\[ \Vert u(t,\cdot)\Vert \leq C_{1}(\Vert u_{0}\Vert + \Vert u_{1}\Vert_{\infty}) + C_{2}\sqrt{\log t}\Vert u_{1}\Vert_{1}\quad (t \gg 1),\]
where $C_{j} > 0$ {\rm (}$j = 1,2${\rm )} are some constants depending on $L$.
\end{pro}
\begin{rem}{\rm When $K(x)$ is a constant function, in \cite{JHDE-ike} one has already obtained the estimate $\Vert u(t,\cdot)\Vert \sim \sqrt{\log t}$ as $t \to \infty$, so that this upper bound derived in Proposition \ref{2-7} may be optimal. For variable coefficients, obtaining a lower bound would be quite difficult. Note that the optimality of this result is only assured with certainty due to the results of Paper \cite{JHDE-ike}.}
\end{rem}

In order to prove Proposition \ref{2-7}, for the smooth solution $u(t,x)$ to problem \eqref{0eqn}-\eqref{0initial}, as in \cite{IM} we set 
$$v(t,x) := \int_{0}^{t}u(s,x)ds.$$
Then, the function $v(t,x)$ satisfies
\begin{align}
& v_{tt} - \nabla \cdot(K(x)\nabla v) = u_{1},\ \ \ (t,x)\in (0,\infty)\times {\bf R}^{n},\label{6eqn}\\
& v(0,x)= 0, \quad  v_{t}(0,x)= u_{0}(x),\ x\in{\bf R}^{n},\label{6initial}
\end{align}
and 
\[{\rm supp}\,v(t,\cdot) \subset B_{L+k_{1}t} \quad (\forall t \geq 0).\]
Now, we set 
\begin{equation}\label{N1}
h(x) := -\frac{1}{2\pi}\int_{{\bf R}^{2}}\log(\vert x- y\vert)u_{1}(y)dy.
\end{equation}
Then we see that $h \in C^{2}({\bf R}^{2})$, and the function $h(x)$ satisfies the Poisson equation:
\begin{equation}\label{N2}
-\Delta h(x) = u_{1}(x), \quad x \in {\bf R}^{2}. 
\end{equation}  
For these facts, see e.g., \cite[page 23, Theorem 1]{E}. 

We first prepare the following lemma. This lemma has already been stated in \cite{RI}, but it is, of course, also a well-known fact.
\begin{lem}\label{N3}
The function $h(x)$ defined in \eqref{N1} satisfies
\[\vert x\vert \vert \nabla h(x)\vert \leq C\Vert u_{1}\Vert_{L^{1}}\]
for $\vert x \vert \geq 2L$, where $C > 0$ is a constant.
\end{lem}
{\it Proof of Lemma \ref{N3}.}\,\,First of all, note that under the assumption $\vert x\vert \geq 2L$ we can get
\[\vert x-y\vert \geq \vert x\vert - L \geq L,\quad \vert x\vert -L \geq \frac{1}{2}\vert x\vert\]
for $y \in {\bf R}^{2}$ satisfying $\vert y\vert \leq L$. After simple elementary computations on \eqref{N1} and the compact support assumption on the initial datum $u_{1}$ one can get the estimate
\[\vert\nabla h(x)\vert \leq \frac{1}{2\pi}\int_{{\bf R}^{2}}\frac{\vert u_{1}(y)\vert}{\vert x-y\vert}dy\]
\[= \frac{1}{2\pi}\int_{\vert y\vert \leq L}\frac{\vert u_{1}(y)\vert}{\vert x-y\vert}dy \leq \frac{1}{2\pi}\frac{1}{\vert x\vert-L}\int_{{\bf R}^{2}}\vert u_{1}(y)\vert dy \leq \frac{1}{\pi}\frac{1}{\vert x\vert}\Vert u_{1}\Vert_{L^{1}},\]
which implies the desired estimate with a constant $C := \frac{1}{\pi}$.
\hfill
$\Box$

The following lemma is essentially the same as what was already presented in \cite{RI}, with only minor adjustments made to fit the given equation.
\begin{lem}\label{N10} For the function $h(x)$ defined in \eqref{N1} it holds that
\[\int_{\vert x\vert \leq  2L + k_{1}t}\vert\nabla h(x)\vert^{2} dx \leq I_{h} + 2\pi C^{2}\Vert u_{1}\Vert_{L^{1}}^{2} \log (2L+k_{1}t) \quad (t \geq 0),\]
where
\[I_{h} := \int_{\vert x\vert \leq 2L}\vert\nabla h(x)\vert^{2}dx.\]
\end{lem}
{\it Proof.}\,By multiplying both sides of \eqref{6eqn} by $v_{t}$, and integrating over $[0,t]\times {\bf R}^{2}$ it follows that
\[\frac{1}{2}\Vert v_{t}(t,\cdot)\Vert_{L^{2}}^{2} + \frac{1}{2}\Vert\sqrt{K(\cdot)}\nabla v(t,\cdot)\Vert_{L^{2}}^{2}\]
\begin{equation}\label{N4}
= \frac{1}{2}\Vert u_{0}\Vert_{L^{2}}^{2} + \int_{{\bf R}^{2}}u_{1}(x)v(t,x)dx.
\end{equation}
Then, from the finite speed of propagation property for $v(t,x)$, \eqref{N2} and the integration by parts we see that
\[\left\vert \int_{{\bf R}^{2}}u_{1}(x)v(t,x)dx\right\vert = \left\vert \int_{\vert x\vert \leq 2L + k_{1}t}u_{1}(x)v(t,x)dx \right\vert\]
\[= \left\vert -\int_{\vert x\vert \leq 2L + k_{1}t}\Delta h(x)v(t,x)dx \right\vert = \left\vert \int_{\vert x\vert \leq 2L + k_{1}t}\nabla h(x)\cdot\nabla v(t,x)dx \right\vert\]
\[\leq \int_{\vert x\vert \leq 2L + k_{1}t}\vert\nabla h(x)\vert\vert\nabla v(t,x)\vert dx\]
\[\leq C_{\varepsilon}\int_{\vert x\vert \leq 2L + k_{1}t}\vert\nabla h(x)\vert^{2}dx + \varepsilon\int_{\vert x\vert \leq 2L + k_{1}t}\vert\nabla v(t,x)\vert^{2}dx\]
\begin{equation}\label{N5}
= C_{\varepsilon}\int_{\vert x\vert \leq 2L + k_{1}t}\vert\nabla h(x)\vert^{2}dx + \varepsilon\int_{{\bf R}^{2}}\vert\nabla v(t,x)\vert^{2}dx
\end{equation} 
with some parameter $\varepsilon > 0$ and a constant $C_{\varepsilon} > 0$ dependng only on $\varepsilon > 0$, where one has just used the facts that $v(t,x) = 0$, and $\vert\nabla h(x)\vert < +\infty$ for $\vert x\vert = 2L + k_{1}t$ for the boundary integral.
Thus, {\bf (K-1)}, \eqref{N4} and \eqref{N5} imply the estimate
\[\frac{1}{2}\Vert v_{t}(t,\cdot)\Vert_{L^{2}}^{2} + (\frac{k_{m}}{2}-\varepsilon)\Vert\nabla v(t,\cdot)\Vert_{L^{2}}^{2}\]
\begin{equation}\label{N6}
\leq \frac{1}{2}\Vert u_{0}\Vert_{L^{2}}^{2} + C_{\varepsilon}\int_{\vert x\vert \leq 2L + k_{1}t}\vert\nabla h(x)\vert^{2}dx. 
\end{equation}
Note that 
\[\int_{\vert x\vert \leq 2L + k_{1}t}\vert\nabla h(x)\vert^{2}dx < +\infty\]
for each $t \geq 0$ because of $h \in C^{2}({\bf R}^{2})$. By choosing the parameter $\varepsilon > 0$ small enough one can get the crucial estimate because of $v_{t} = u$.
\begin{equation}\label{N7}
\Vert u(t,\cdot)\Vert_{L^{2}}^{2} \leq \Vert u_{0}\Vert_{L^{2}}^{2} + 2C_{\varepsilon}\int_{\vert x\vert \leq 2L + k_{1}t}\vert\nabla h(x)\vert^{2}dx.
\end{equation}
Let us estimate the second term of the right hand side of \eqref{N7} by using Lemma \ref{N3}. 
First, we see that
\begin{equation}\label{N8}
\int_{\vert x\vert \leq 2L + k_{1}t}\vert\nabla h(x)\vert^{2}dx = I_{h} + \int_{2L \leq \vert x\vert \leq 2L + k_{1}t}\vert\nabla h(x)\vert^{2}dx.
\end{equation}
Thus, from Lemma \ref{N3} it follows that
\[\int_{2L \leq \vert x\vert \leq 2L + k_{1}t}\vert\nabla h(x)\vert^{2}dx \leq C^{2}\Vert u_{1}\Vert_{L^{1}}^{2}\int_{2L \leq \vert x\vert \leq 2L + k_{1}t}\frac{1}{\vert x\vert^{2}}dx\]
\[= 2\pi C^{2}\Vert u_{1}\Vert_{L^{1}}^{2}\int_{2L}^{2L + k_{1}t}\frac{r}{r^{2}}dr = 2\pi C^{2}\Vert u_{1}\Vert_{L^{1}}^{2}\left(\log(2L+k_{1}t)-\log(2L)\right)\]
\begin{equation}\label{N9}
\leq 2\pi C^{2}\Vert u_{1}\Vert_{L^{1}}^{2}\log(2L+k_{1}t)\quad (t \geq 0). 
\end{equation}
Therefore, by \eqref{N8} and \eqref{N9} one has arrived at the desired estimate.
\hfill
$\Box$

Finally, let us estimate $I_{h}$ in terms of the initial velocity $u_{1}$. Note that generally it holds that
\[\int_{{\bf R}^{2}}\vert\nabla h(x)\vert^{2}dx = +\infty.\] 
It is important to point out that the integral region of the quantity $I_{h}$ is localized to $\vert x\vert \leq 2L$ as is well-known. Of course, we know that its integral value is finite because of $h \in C^{2}({\bf R}^{2})$, but we will try to suppress it by a certain quantity of initial velocity from above.\\
The lemma below has also already been discussed in \cite{RI}, but it is included here for the reader's convenience.
\begin{lem}\label{N13}\, For the function $h(x)$ defined in \eqref{N1} it holds that 
\[I_{h} = \int_{\vert x\vert \leq 2L}\vert\nabla h(x)\vert^{2}dx \leq C_{L}\Vert u_{1}\Vert_{L^{\infty}}^{2}\]
with some constant $C_{L} > 0$ depending only on $L$.
\end{lem}
{\it Proof.}\,For an arbitrarily fixed $x_{0} \in B_{2L}$, we can easily arrive at the intermediate estimate, where $B_{r}(a) := \{x \in {\bf R}^{2}\,:\,\vert x-a\vert \leq r\}$, and we set $B_{r}(0) = B_{r}$.
\[\vert\nabla h(x_{0})\vert \leq \frac{1}{2\pi}\int_{{\bf R}^{2}}\frac{\vert u_{1}(y)\vert}{\vert x_{0}-y\vert}dy = \frac{1}{2\pi}\int_{\vert y\vert \leq L}\frac{\vert u_{1}(y)\vert}{\vert x_{0}-y\vert}dy \]
\[= \frac{1}{2\pi}\int_{B_{L}(x_{0})}\frac{\vert u_{1}(x_{0}-z)\vert}{\vert z\vert}dz\]
\begin{equation}\label{N11}
\leq \frac{\Vert u_{1}\Vert_{L^{\infty}}}{2\pi}\int_{B_{L}(x_{0})}\frac{1}{\vert z\vert}dz.
\end{equation}
It should be noticed that in the case of $x_{0} \in B_{2L}$, we have $B_{L}(x_{0}) \subset B_{4L}$ sufficiently. Thus, we get
\begin{equation}\label{N12}
\int_{B_{L}(x_{0})}\frac{1}{\vert z\vert}dz \leq \int_{B_{4L}}\frac{1}{\vert z\vert}dz = 8\pi L.
\end{equation}
Put $x_{0}$ back into $x$ to get the following by \eqref{N11} and \eqref{N12}. 
\[\vert\nabla h(x)\vert \leq 4L\Vert u_{1}\Vert_{L^{\infty}}\quad \forall x \in B_{2L}.\]
Therefore, one can obtain the desired estimate.
\hfill
$\Box$

Under these preparations, let us prove Proposition \ref{2-7}.\\ 

{\it Proof of Proposition \ref{2-7}.}\,\,Indeed, the desired estimate is a direct consequence of \eqref{N7} and Lemmas \ref{N10} and \ref{N13}.
\hfill
$\Box$
\begin{rem}{\rm The key insight in the proof of Proposition \ref{2-7} lies in recognizing that we solve 
\[-\Delta h(x) = u_{1}(x)\]
rather than 
\begin{equation}\label{N14}
-\nabla\cdot(K(x)\nabla h(x)) = u_{1}(x)
\end{equation}
like Morawetz's original idea (see \cite{Mora}). Solving \eqref{N14} generally seems difficult because the form of $K(x)$ is unknown, making it hard to estimate solutions. As can be seen from the above proof, it should be noted that the mere "existence" of a solution to \eqref{N14} does not lead to the desired resolution of the problem. Furthermore, removing the compactness condition on the initial value support from the proof appears challenging.}
\end{rem}

\section{Proof of Theorem \ref{theorem}}

Under several preparations in Section $2$, let us prove Theorem \ref{theorem} in this section.

First of all, by applying the assumption {\bf (K-2)} to the identity derived in Lemma \ref{2-1} one can get the inequality:
\begin{equation}\label{3-1}
tE_{u}(t) \leq J_{0} -\frac{1}{2}\int_{{\bf R}^{2}}u_{t}(t,x)u(t,x)dx - \int_{{\bf R}^{2}}u_{t}(t,x)(x\cdot \nabla u(t,x))dx
\end{equation}
for $t \geq 0$. 

We will now evaluate the second and third terms on the right-hand side of the inequality \eqref{3-1} in the required form. Tobegin with, we get the following lemma.
\begin{lem}\label{3-2}\, The smooth solution $u(t,x)$ to problem \eqref{0eqn}-\eqref{0initial} satisfies the following growth property.
\[\int_{{\bf R}^{2}}\vert u_{t}(t,x)u(t,x)\vert dx \leq C_{3}\sqrt{E_{u}(0)}\left( \Vert u_{0}\Vert + \Vert u_{1}\Vert_{\infty} + \Vert u_{1}\Vert_{1} \right)\sqrt{\log t}\]
for $t \gg 1$, where the constant $C_{3} > 0$ depends on $L > 0$, and does not depend on any norms of initial data.
\end{lem}
{\it Proof.} Since $\Vert u_{t}(t,\cdot)\Vert \leq \sqrt{2E_{u}(t)} = \sqrt{2E_{u}(0)}$, by the Cauchy-Schwarz inequality it follows that
\[\int_{{\bf R}^{2}}\vert u_{t}(t,x)u(t,x)\vert dx \leq \Vert u(t,\cdot)\Vert\Vert u_{t}(t,\cdot)\Vert \leq \sqrt{2E_{u}(0)}\Vert u(t,\cdot)\Vert.\]
By applying the result obtained in Proposition \ref{2-7} it holds that
\[\int_{{\bf R}^{2}}\vert u_{t}(t,x)u(t,x)\vert dx \leq \sqrt{2E_{u}(0)}\{C_{1}(\Vert u_{0}\Vert + \Vert u_{1}\Vert_{\infty}) + C_{2}\sqrt{\log t}\Vert u_{1}\Vert_{1}\}\]
for $t \gg 1$. This implies the desired estimate.
\hfill
$\Box$

Next, let us evaluate the third term on the right-hand side of \eqref{3-1} based on Lemma \ref{2-6}. We transform using the form of the weight function $\psi(t,x)$ defined in Section 2. 
\begin{lem}\label{3-3}\, Let $R > r_{0}$, and $t > \frac{R}{\sqrt{k_{0}}}$. Then, the smooth solution $u(t,x)$ to problem \eqref{0eqn}-\eqref{0initial} satisfies
\[\int_{{\bf R}^{2}}\left\vert u_{t}(t,x)(x\cdot \nabla u(t,x))\right\vert dx \leq \frac{R}{\sqrt{k_{m}}}E_{u,R}(t) + \frac{C_{r_{0}}}{\sqrt{k_{m}}}\int_{{\bf R}^{2}}(1+\vert x\vert)e(0,x)dx + t\int_{\vert x\vert \geq R}e(t,x)dx.\]
\end{lem}
{\it Proof.}\,In fact, we see that
\[\int_{{\bf R}^{2}}\left\vert u_{t}(t,x)(x\cdot \nabla u(t,x))\right\vert dx\]
\[ = \frac{1}{\sqrt{k_{m}}}\int_{B_{R}}\sqrt{k_{m}}\vert u_{t}(t,x)\vert\vert x\vert\vert\nabla u(t,x)\vert dx + \int_{\vert x\vert \geq R}\vert u_{t}(t,x)\vert\vert x\vert\vert\nabla u(t,x)\vert dx\]
\[\leq \frac{R}{\sqrt{k_{m}}}\int_{B_{R}}e(t,x)dx + \left(\int_{\vert x\vert \geq \sqrt{k_{0}}t} + \int_{\sqrt{k_{0}}t \geq \vert x\vert \geq R}\right)\vert u_{t}(t,x)\vert\vert x\vert\vert\nabla u(t,x)\vert dx\]
\[= \frac{R}{\sqrt{k_{m}}}E_{u,R}(t) + \int_{\vert x\vert \geq \sqrt{k_{0}}t}\vert u_{t}(t,x)\vert(\vert x\vert - \sqrt{k_{0}}t)\vert\nabla u(t,x)\vert dx\]
\[+ \sqrt{k_{0}}t\int_{\vert x\vert \geq \sqrt{k_{0}}t}\vert u_{t}(t,x)\vert \vert\nabla u(t,x)\vert dx + \sqrt{k_{0}}t\int_{\sqrt{k_{0}}t\geq \vert x\vert \geq R}\vert u_{t}(t,x)\vert\vert\nabla u(t,x)\vert dx\]
\[\leq \frac{R}{\sqrt{k_{m}}}E_{u,R}(t) + \frac{1}{\sqrt{k_{m}}}\int_{\vert x\vert \geq \sqrt{k_{0}}t}\psi(t,x)\vert u_{t}(t,x)\vert\sqrt{k_{m}}\vert\nabla u(t,x)\vert dx + t\int_{\vert x\vert \geq R}\vert u_{t}(t,x)\vert(\sqrt{k_{0}}\vert\nabla u(t,x)\vert) dx\]
\[\leq \frac{R}{\sqrt{k_{m}}}E_{u,R}(t) + \frac{1}{\sqrt{k_{m}}}\int_{\vert x\vert \geq \sqrt{k_{0}}t}\psi(t,x)e(t,x)dx + t\int_{\vert x\vert \geq R}e(t,x)dx.\]
The estimate of the last two lines essentially uses condition {\bf (K-3)}. This implies the desired estimate by using Lemma \ref{2-6} and \eqref{2-3-2}.
\hfil
$\Box$
\begin{rem}{\rm The discussion appeared in the proof of Lemma \ref{3-2} is independent of the size $L$ of support of the initial data. This may serve as a hint when generalizing the initial value class.}
\end{rem}

Let us finalize the proof of Theorem \ref{theorem}.\\

{\it Proof of Theorem \ref{theorem}}\, Let $R > r_{0}$ and $t \gg 1 $ and $t >\frac{R}{\sqrt{k_{0}}}$. Then, it follows from Lemmas \ref{3-2} and \ref{3-3} and \eqref{3-1} it follows that
\[tE_{u,R}(t) + t\int_{\vert x\vert \geq R}e(t,x)dx = tE_{u}(t)\] 
\[\leq J_{0} + C_{3}\sqrt{E_{u}(0)}\left( \Vert u_{0}\Vert + \Vert u_{1}\Vert_{\infty} + \Vert u_{1}\Vert_{1} \right)\sqrt{\log t}\]
\[+ \frac{R}{\sqrt{k_{m}}}E_{u,R}(t) + \frac{C_{r_{0}}}{\sqrt{k_{m}}}\int_{{\bf R}^{2}}(1+\vert x\vert)e(0,x)dx + t\int_{\vert x\vert \geq R}e(t,x)dx,\]
which implies
\[(t - \frac{R}{\sqrt{k_{m}}})E_{u,R}(t) \leq J_{0} + C_{3}\sqrt{E_{u}(0)}\left( \Vert u_{0}\Vert + \Vert u_{1}\Vert_{\infty} + \Vert u_{1}\Vert_{1} \right)\sqrt{\log t}\]
\[+ \frac{C_{r_{0}}}{\sqrt{k_{m}}}\int_{{\bf R}^{2}}(1+\vert x\vert)e(0,x)dx \quad (t \gg 1).\]
This shows the desired statement of Theorem \ref{theorem}.
\hfill
$\Box$

\section{Generalization of the Conditions for $K(x)$}

In this section, we generalize the condition {\bf (K-2)} assuming monotonicity. Therefore, we will replace {\bf (K-2)} with the following condition.\\
{\bf (K-4)}\, There exists a real number $\gamma_{0} \in [0,1)$ such that $x\cdot\nabla K(x) \leq \gamma_{0}K(x)$ for all $x \in {\bf R}^{2}$.\\
One can state the following slightly generalized version of Theorem \ref{theorem}.
\begin{theo}\label{theorem2}\,Assume \eqref{I-1} on the initial data $u_{j} \in C_{0}^{\infty}({\bf R}^{2})$ {\rm (}$j = 0,1${\rm )} and the conditions {\bf (K-1)}, {\bf (K-3)} and {\bf (K-4)}. Let $R > r_{0}$ be an arbitrarily fixed real number. Then, for the solution $u(t,x)$ to problem \eqref{0eqn}-\eqref{0initial} it holds that
\[E_{u,R}(t) = O(t^{\gamma_{0}-1}\sqrt{\log t}) \quad (t \gg 1).\] 
\end{theo}
{\it Proof.}\, Let $R > r_{0}$. Then, as in the proof of Theorem \ref{theorem} it follows from {\bf (K-3)} and Lemma \ref{2-1} that
\[tE_{u}(t) \leq J_{0} + C_{3}\sqrt{E_{u}(0)}\left( \Vert u_{0}\Vert + \Vert u_{1}\Vert_{\infty} + \Vert u_{1}\Vert_{1} \right)\sqrt{\log t}\]
\[+ \frac{R}{\sqrt{k_{m}}}E_{u,R}(t) + \frac{C_{r_{0}}}{\sqrt{k_{m}}}\int_{{\bf R}^{2}}(1+\vert x\vert)e(0,x)dx + t\int_{\vert x\vert \geq R}e(t,x)dx + \gamma_{0}\int_{0}^{t}E_{u,R}(s)ds.\]
This implies
\begin{equation}\label{4-1}
(t - \frac{R}{\sqrt{k_{m}}})E_{u,R}(t) \leq A_{0}(t) + \gamma_{0}\int_{0}^{t}E_{u,R}(s)ds \quad (t \gg 1),
\end{equation}
where
\[A_{0}(t) := J_{0} + C_{3}\sqrt{E_{u}(0)}\left( \Vert u_{0}\Vert + \Vert u_{1}\Vert_{\infty} + \Vert u_{1}\Vert_{1} \right)\sqrt{\log t} + \frac{C_{r_{0}}}{\sqrt{k_{m}}}\int_{{\bf R}^{2}}(1+\vert x\vert)e(0,x)dx.\]
Let us solve the differential inequality \eqref{4-1}. For this set
\[w(t) := (t-\frac{R}{\sqrt{k_{m}}})^{-\gamma_{0}}\int_{0}^{t}E_{u,R}(s)ds.\]
Then, by differentiating $w(t)$ one finds that
\[w'(t) = (t-\frac{R}{\sqrt{k_{m}}})^{-\gamma_{0}-1}\left((t-\frac{R}{\sqrt{k_{m}}})E_{u,R}(t) - \gamma_{0}\int_{0}^{t}E_{u,R}(s)ds\right)\] 
\[\leq A_{0}(t)(t-\frac{R}{\sqrt{k_{m}}})^{-\gamma_{0}-1}\quad (a.e.\,\, t \gg 1).\]
Integration it over $[t_{0},t]$ with large fixed $t_{0} \gg 1$, we see that
\[w(t) \leq w(t_{0}) + \left(J_{0} + \frac{C_{r_{0}}}{\sqrt{k_{m}}}\int_{{\bf R}^{2}}(1+\vert x\vert)e(0,x)dx\right)\int_{t_{0}}^{t}(s-\frac{R}{\sqrt{k_{m}}})^{-\gamma_{0}-1}ds\]
\[+ \left(C_{3}\sqrt{E_{u}(0)}( \Vert u_{0}\Vert + \Vert u_{1}\Vert_{\infty} + \Vert u_{1}\Vert_{1})\right)\int_{t_{0}}^{t}\sqrt{\log s}(s-\frac{R}{\sqrt{k_{m}}})^{-\gamma_{0}-1}ds \quad (t \geq t_{0}).\]
Here, note that
\[\int_{t_{0}}^{t}\sqrt{\log s}(s-\frac{R}{\sqrt{k_{m}}})^{-\gamma_{0}-1}ds \leq \sqrt{\log t}\int_{t_{0}}^{t}(s-\frac{R}{\sqrt{k_{m}}})^{-\gamma_{0}-1}ds\]
\[\leq \frac{1}{\gamma_{0}}\sqrt{\log t}(t_{0}-\frac{R}{\sqrt{k_{m}}})^{-\gamma_{0}} \quad (t \geq t_{0}).\]
Therefore, one can get
\[w(t) \leq w(t_{0}) + \left(J_{0} + \frac{C_{r_{0}}}{\sqrt{k_{m}}}\int_{{\bf R}^{2}}(1+\vert x\vert)e(0,x)dx\right)\frac{1}{\gamma_{0}}(t_{0}-\frac{R}{\sqrt{k_{m}}})^{-\gamma_{0}}\]
\begin{equation}\label{4-2}
+  \left(C_{3}\sqrt{E_{u}(0)}( \Vert u_{0}\Vert + \Vert u_{1}\Vert_{\infty} + \Vert u_{1}\Vert_{1})\right)\frac{1}{\gamma_{0}}(t_{0}-\frac{R}{\sqrt{k_{m}}})^{-\gamma_{0}}\sqrt{\log t}.
\end{equation}
Additionally, from \eqref{4-1} we can get
\[(t-\frac{R}{\sqrt{k_{m}}})E_{u,R}(t) \leq A_{0}(t) + \gamma_{0}\int_{0}^{t}E_{u,R}(s)ds\]
\begin{equation}\label{4-3}
= A_{0}(t) + \gamma_{0}(t-\frac{R}{\sqrt{k_{m}}})^{\gamma_{0}}w(t).  
\end{equation}
By combining \eqref{4-2} and \eqref{4-3} it follows that
\[(t-\frac{R}{\sqrt{k_{m}}})E_{u,R}(t) \leq A_{0}(t) + \gamma_{0}(t-\frac{R}{\sqrt{k_{m}}})^{\gamma_{0}}w(t_{0})\]
\[+ \left(J_{0} + \frac{C_{r_{0}}}{\sqrt{k_{m}}}\int_{{\bf R}^{2}}(1+\vert x\vert)e(0,x)dx\right)(t_{0}-\frac{R}{\sqrt{k_{m}}})^{-\gamma_{0}}(t-\frac{R}{\sqrt{k_{m}}})^{\gamma_{0}}\]
\[+ \left(C_{3}\sqrt{E_{u}(0)}( \Vert u_{0}\Vert + \Vert u_{1}\Vert_{\infty} + \Vert u_{1}\Vert_{1})\right)(t_{0}-\frac{R}{\sqrt{k_{m}}})^{-\gamma_{0}}(t-\frac{R}{\sqrt{k_{m}}})^{\gamma_{0}}\sqrt{\log t}.\]
This gives 
\[E_{u,R}(t) \leq A_{0}(t)(t-\frac{R}{\sqrt{k_{m}}})^{-1} + \gamma_{0}(t-\frac{R}{\sqrt{k_{m}}})^{\gamma_{0}-1}w(t_{0})\]
\[+ \left(J_{0} + \frac{C_{r_{0}}}{\sqrt{k_{m}}}\int_{{\bf R}^{2}}(1+\vert x\vert)e(0,x)dx\right)(t_{0}-\frac{R}{\sqrt{k_{m}}})^{-\gamma_{0}}(t-\frac{R}{\sqrt{k_{m}}})^{\gamma_{0}-1}\]
\[+ \left(C_{3}\sqrt{E_{u}(0)}( \Vert u_{0}\Vert + \Vert u_{1}\Vert_{\infty} + \Vert u_{1}\Vert_{1})\right)(t_{0}-\frac{R}{\sqrt{k_{m}}})^{-\gamma_{0}}(t-\frac{R}{\sqrt{k_{m}}})^{\gamma_{0}-1}\sqrt{\log t},\]
which implies the desired estimate
\[E_{u,R}(t) = O(t^{\gamma_{0}-1}\sqrt{\log t})\quad (t \gg 1).\]
Note that $A_{0}(t) \sim \sqrt{\log t}$. 
\hfill
$\Box$
\begin{rem}{\rm To justify the proof of the above theorem, note the following.

Since the class of solutions under consideration is at least $u \in C([0,\infty);H^{1}({\bf R}^{2}))\cap C^{1}([0,\infty);L^{2}({\bf R}^{2}))$, the continuity of the function $t \mapsto \Vert u_{t}(t)\Vert_{L^{2}(B_{R})}$ for each $R > 0$ follows from the following inequality. For fixed $t_{1} \geq 0$
\[\left\vert \Vert u_{t}(t)\Vert_{L^{2}(B_{R})}-\Vert u_{t}(t_{1})\Vert_{L^{2}(B_{R})} \right\vert \leq \Vert u_{t}(t)-u_{t}(t_{1})\Vert_{L^{2}({\bf R}^{2})}.\]
Similarly, the continuity of the function $t \mapsto \Vert\sqrt{K(\cdot)}\nabla u(t)\Vert_{L^{2}(B_{R})}$ for each $R > 0$ can be checked under the condition that $K \in L^{\infty}({\bf R}^{2})$. These imply the continuity of the mapping $t \mapsto E_{u,R}(t)$.}
\end{rem}
\begin{rem}{\rm One can give one simple example for $K \in {\rm BC}^{1}({\bf R}^{2})$ satisfying the conditions {\bf (K-1)}, {\bf (K-3)} and {\bf (K-4)}:\\
\[ K(x) = \left\{
     \begin{array}{ll}
       \displaystyle{(1+\vert x\vert^{2})^{\frac{\gamma_{0}}{2}}}&
              \qquad \vert x\vert \leq r_{0}/2, \\[0.2cm]
       \displaystyle{(1 + r_{0}^{2})^{\frac{\gamma_{0}}{2}}}& \qquad  r_{0} < \vert x\vert,
\end{array} \right. 
\]
where $\gamma_{0} \in [0,1)$ is a given parameter. It suffices to connect the graph of $K(x)$ smoothly and gradually in the region $r_{0}/2 \leq \vert x\vert \leq r_{0}$. A slight increase in $K(x)$ is also permissible near the origin.}
\end{rem}


\section{Lipschitz perturbation case for $K(x)$}

In this section, we proceed with the same considerations as before when $K(x)$ is Lipschitz continuous. Note that Lemma \ref{2-1} continues to hold. 

First, we assume the following for $K(x)$ newly.\\
{\bf (K-5)}\, $\nabla K \in (L^{\infty}({\bf R}^{2}))^{2}$.\\
Note that the condition {\bf (K-5)} means that $K(x)$ is Lipschitz continuous on ${\bf R}^{2}$.

One can state the following one more generalized version of Theorem \ref{theorem}.
\begin{theo}\label{theorem3}\,Assume \eqref{I-1} on the initial data $u_{j} \in C_{0}^{\infty}({\bf R}^{2})$ {\rm (}$j = 0,1${\rm )} and the conditions {\bf (K-1)}, {\bf (K-3)} and {\bf (K-5)}. Let $R > r_{0}$ be an arbitrarily fixed real number. Then, for the solution $u(t,x)$ to problem \eqref{0eqn}-\eqref{0initial} it holds that
\[E_{u,R}(t) = O(t^{\eta_{0}-1}\sqrt{\log t}) \quad (t \gg 1),\]
provided that 
\[\eta_{0} := \frac{r_{0}}{k_{m}}\Vert\nabla K\Vert_{\infty} \in [0,1).\] 
\end{theo}
{\it Proof.} Based on Lemma \ref{2-1}, the following evaluation should be noted; the rest can be computed similarly to the previous section. In fact, using {\bf (K-1)}, {\bf (K-3)} and {\bf (K-5)} we see that
\[\frac{1}{2}\int_{0}^{t}\int_{{\bf R}^{2}}\vert (x\cdot\nabla K(x))\vert\vert\nabla u(s,x)\vert^{2}dxds = \frac{1}{2}\int_{0}^{t}\int_{B_{r_{0}}}\vert (x\cdot\nabla K(x))\vert\vert\nabla u(s,x)\vert^{2}dxds\]
\[\leq \frac{r_{0}}{2k_{m}}\Vert \nabla K\Vert_{\infty}\int_{0}^{t}\int_{B_{r_{0}}}k_{m}\vert\nabla u(s,x)\vert^{2}dxds \leq \frac{r_{0}}{k_{m}}\Vert \nabla K\Vert_{\infty}\int_{0}^{t}E_{u,R}(s)ds.\]
Therefore, following the same procedure as in the previous Section 4, we obtain the following.
\begin{equation}\label{5-1}
(t - \frac{R}{\sqrt{k_{m}}})E_{u,R}(t) \leq A_{0}(t) + \eta_{0}\int_{0}^{t}E_{u,R}(s)ds \quad (t \gg 1),
\end{equation}
which implies the desired estimate (see \eqref{4-1} and the subsequent calculations).
\hfill
$\Box$
\begin{rem}{\rm Note that although we assumed $K \in {\rm BC}^{1}({\bf R}^{2})$ for simplicity regarding $K(x)$ at the beginning of the paper, $K \in W^{1,\infty}({\bf R}^{2})$ is actually sufficient. As the results show, even in two-dimensional space, imposing a certain constraint on the magnitude of the derivative of the variable coefficient $K(x)$ allows us to derive that the local energy decays at an algebraic order. Although the equations differ slightly, a comparison with Paper \cite{Sha} should be interesting. Further calculations could be used to develop a similar theory for the following equation with new density function $\rho(x) > 0$
\[\rho(x)^{2}u_{tt}(t,x)-\nabla \cdot(K(x)\nabla u(t,x)) = 0,\]
but this is left to the reader's interest (see \cite{A} as a hint).}
\end{rem}






\end{document}